\documentclass[11pt]{amsart}
\usepackage{amsmath, amssymb, amsfonts, amscd,mathptmx}
\usepackage{epsfig}

\newcommand{\n}{\mbox{$\mathfrak{n}$}}
\newcommand{\lf}{\mbox{$\mathfrak{l}$}}

\newcommand{\m}{\mbox{$\mathfrak{m}$}}

\newtheorem{theorem}{Theorem}[section]

\newtheorem{lemma}[theorem]{Lemma}
\newtheorem{remark}[theorem]{Remark}

\begin{document}

\title[nilpotent Lie algebras  associated with graphs]{graphs and two-step Nilpotent Lie algebras}

\author{Meera G. Mainkar}

\address{Department of Mathematics, Pearce Hall, Central Michigan University, Mt. Pleasant, MI 48859, USA} \email{maink1m@cmich.edu}

\thanks{
  The author was funded by Early Career Investigator Grant C61940 from Central Michigan University.}

\maketitle
\begin{abstract}

 We consider a  method popular in the literature of associating a two-step nilpotent Lie algebra with a finite simple graph.
  We prove that the two-step  nilpotent Lie algebras  associated with two graphs
    are Lie isomorphic if and only
  if the graphs from which they arise are isomorphic.
\end{abstract}

\section{Introduction}\label{intro}
Recently in connection with the study of some very interesting geometrical and dynamical properties of nilmanifolds,
 the simply connected nilpotent Lie groups or equivalently  their Lie algebras  have received the attention of
many researchers.  Among these  nilpotent Lie algebras, the two-step ones are the simplest (after the abelian ones) and most widely studied \cite{D-M, DDM, De, D-D, D-V1, Dem, E1, E2, F1, F2, GorM, GM, L-W1, L-W2, P, P-T}.

 Let $k$ be a  field with characteristic  not equal to 2. We recall that
a  non-abelian finite dimensional Lie algebra $\n$ over $k$ is said to be a {\em two-step nilpotent Lie algebra} if $[\n, [\n, \n]] = \{0\}$.
 Every two-step nilpotent Lie algebra $\n$ over $k$ can be realized as a vector space direct sum $V \oplus (\bigwedge^2V)/W$ where $V$ is a finite dimensional $k$-vector space
  and $W$ is a subspace of the exterior power $\bigwedge^2V$.   The  Lie bracket structure on $\n$ is  given by:

  \begin{enumerate}

  \item $[v_1, v_2] = v_1 \wedge v_2 \text { mod } W \,\,\text{ for } v_1, v_2 \in V,$ and

 \item  $[x, y] = 0 \text{ for } x \in \n \text{ and } y \in (\bigwedge^2V)/W$
  \end{enumerate}

A combinatorial approach for construction of  two-step nilpotent Lie algebras was described in \cite{D-M}.
 Subsequently this construction has been used by many authors \cite{F1,F2,P-T,L-W2}. We recall the construction.
 Let $(S, E)$ be a finite simple graph where $S$ is the set of vertices and $E$ is the set of edges, where it is assumed that
 there are no loops and no multiple edges connecting the same pair of vertices. We associate with $(S, E)$ a two-step nilpotent Lie algebra $\n = \n(S, E)$ over $k$ in the following way.
  The underlying vector space of $\n$ is $V \oplus  (\bigwedge^2V)/W$ where $V$ is the $k$-vector space consisting of formal $k$-linear combinations of elements of $S$ (so that $S$ is a basis of $V$), and $W$ is the subspace of the exterior power $\bigwedge^2V$ spanned by
   the vectors $\alpha \wedge \beta$ where $\alpha \beta \notin E$. The Lie bracket structure on $\n$ is defined as above. We note that the space $(\bigwedge^2V)/W$ has dimension $|E|$. In fact it is spanned by $\{\alpha \wedge \beta \text{ mod }  W \mid \alpha \beta \in E\}$.
 Consequently $\n(S, E)$ has dimension $|S| + |E|$.

Note that the construction of $\n(S,E)$ form the graph $(S, E)$ is functorial in the sense that if $f : (S, E) \to (S', E')$ is an isomorphism  of graphs, then we obtain
an isomorphism $f_{*}: \n(S, E) \to \n(S', E')$ in the natural way: $f_*(\alpha) = f(\alpha)$ for each $\alpha \in S$ and for each edge $\alpha \beta \in E$, $f_*(\alpha \wedge \beta) = f(\alpha) \wedge f(\beta)$. It is easy to see that $f_*$ extends linearly to a Lie algebra isomorphism from $\n(S,E)$ to $\n(S', E')$.

In this note, we consider the converse question and prove the following:

  \begin{theorem}\label{main}
  Let $(S, E)$ and $(S', E')$ be finite simple graphs.
  If the two Lie algebras $\n(S, E)$ and $\n(S', E')$ are isomorphic then the graphs $(S, E)$ and $(S', E')$ are also isomorphic.
   \end{theorem}

   This result has already found use in some recent investigations \cite{L-W2, P-T}.
   In \cite{L-W2}, the authors provide a method to construct  Einstein solvmanifolds by using the Lie algebras $\n(S, E)$. Using our main theorem, these solvmanifolds are isometric if and only if the graphs are isomorphic, which gives us examples  of nonisometric Einstein solvmanifolds.

    In \cite{P-T}, the authors construct symplectic two-step nilpotent Lie algebras associated with graphs.
    In this case our main result implies  that  {\em there are exactly five non-isomorphic two-step nilpotent Lie algebras of dimension six associated with graphs in the above manner.} Indeed, there are five  non-isomorphic graphs $(S, E)$ such that $|S| + |E| = 6 $. Now,
  $\dim \n(S, E)= |S| + |E|$.
  Hence using our Theorem \ref{main}, there are exactly five non-isomorphic two-step nilpotent
 Lie algebras of dimension six associated with graphs.
    This fact was used in
     \cite[Remark 1]{P-T}.

 The group of Lie automorphisms of $\n(S, E)$ was determined in terms of the graph $(S, E)$ in \cite{D-M}. Also, Anosov and ergodic automorphisms on corresponding nilmanifolds were studied.
 In \cite{L-W2}, explicit examples and non-examples
  of  Einstein solvmanifolds were constructed using the Lie algebras $\n(S,E)$ (see also \cite{F1}).  A combinatorial construction of the first and second cohomology
   groups of $\n(S,E)$  was given in \cite{P-T}, and was  used to construct symplectic and contact nilmanifolds.

\vspace{0.1cm}

 %addition
 \begin{remark}\label{groups}
 {\rm Suppose the underlying field $ k = \mathbb R$. Let    $N(S, E)$  denote the unique  simply connected nilpotent Lie group corresponding to the Lie algebra $\n(S, E)$. Then  the Lie group exponential map $\exp : \n(S,E) \to N(S, E)$ is a diffeomorphism (see \cite{R}, p. 6).
  We note that the Baker-Campbell-Hausdorff formula (\cite{V}, p. 119)) gives the  multiplication law in $N(S,E)$  as follows:
  \[\exp(x)\exp(y) = \exp (x + y + \frac{1}{2}[x, y])\]
   for all $x, y \in \n(S,E)$.

   More precisely, we can realize  $N(S,E)$ as $\n(S,E)$ (via the exponential map) with the multiplication defined by
   \[ (v_1, x_1). (v_2, x_2) = (v_1 + v_2, x_1 + x_2 + \frac{1}{2}[v_1, v_2]) \]

   for all $v_1, v_2 \in V$ and $x_1, x_2 \in (\bigwedge^2V)/W$.

   Then our Theorem~\ref{main} implies that the simply connected Lie groups $N(S, E)$ and $N(S',E')$ are Lie isomorphic if and only if the graphs $(S,E)$ and $(S', E')$ are  isomorphic. This can been seen by using the fact that the simply connected Lie groups are Lie isomorphic if and only if their Lie algebras are Lie isomorphic (see \cite{Warner}, p. 101).
   }
   \end{remark}

\vspace{0.1cm}

   \begin{remark}\label{referee}
   {\rm There have been some other constructions  of algebraic structures associated with a simple graph. Some, though not all of these are related with the construction considered here.
   In \cite{KLNR}, the authors consider the $K$-algebra  associated with  a simple graph $(S,E)$ over a field $K$  generated by the set of vertices $S$ and with relations $\alpha \beta = \beta \alpha$ if and only if $\alpha \beta \notin E$. They proved that these $K$-algebras are isomorphic if and only the corresponding graphs are isomorphic.
    This result was further used to prove an analogous result for graph groups in \cite{Dr1}. In \cite{Dr1}, the author considers the group associated with graph $(S, E)$ which is defined as the group generated by the set $S$ and with relations $\alpha \beta = \beta \alpha$ if and only if $\alpha \beta \in E$.   Here we remark that the simply connected nilpotent Lie group $N(S, E)$ (as in Remark \ref{groups}) is not finitely generated.

   Later in \cite{DuK1}, the authors introduce the free partially commutative Lie algebra $\mathfrak{l}(S, E)$ associated with $(S,E)$ which is the quotient of the  free Lie algebra on the set $S$ modulo  the Lie ideal
    generated by $\{ [\alpha, \beta]\mid \alpha\beta \notin E \}$. In fact, the authors study these structures in more generality but we will not go into the details here. Following \cite{Dr1, DuK2, KLNR}, it can be shown that these free partially commutative Lie algebras are Lie isomorphic if and only if the corresponding graphs are isomorphic (see Theorem 1 in \cite{Dr1} for example).
    Here we note that the quotient of $\lf(S,E)$ by the Lie ideal $\lf^3$ generated by the set $\{[x, [y, z]]\mid x, y, z \in \lf(S,E)\}$, is Lie isomorphic to our two-step nilpotent Lie algebra $\n(S,E)$.  The isomorphism result for the free partially commutative Lie algebras can be recaptured by using our Theorem~\ref{main} (if the characteristic of the field $K$ is not equal to 2) as follows: Let $\lf= \lf(S,E)$ and $\lf' = \lf(S',E')$  be Lie isomorphic free partially commutative Lie algebras associated with graphs
    $(S, E)$ and $(S', E')$ respectively. Then the quotient Lie algebras  $\lf / \lf^3$ and $\lf'/\lf'^3$ are Lie isomorphic. That means the two-step nilpotent Lie algebras $\n(S, E)$ and $\n(S', E')$ are isomorphic and our Theorem~\ref{main} implies that the graphs $(S, E)$ and $(S', E')$ are isomorphic.}
\end{remark}

\vspace{0.3cm}

\section{Some Preliminaries}\label{prelim}

       Before we  prove  Theorem~\ref{main}, we need a few facts regarding  the structure of the automorphism group of  $\n = \n(S, E)$.
         We denote by ${\rm Aut}(\n)$ the  group of Lie automorphisms  of $\n$ and $T$ the set of automorphisms $\tau \in {\rm Aut}(\n)$
          such that $\tau(V) = V$ where we recall that $V$ is the  $k$-vector space with basis $S$, and $\n = V \oplus (\bigwedge^2 V)/W$.  We denote by
          $G$ the subgroup of ${\rm GL}(V)$ consisting of the linear automorphisms $\tau|_V$ where  $\tau \in T$.

           Next we see that the group $G$ consists of  $g \in {\rm GL }(V)$ such that the subspace $W \subseteq \bigwedge^2 V$ is invariant
           under the induced  natural action $\wedge^2 g$ of $g$ on $\bigwedge^2V$.
            Indeed,
           if $g \in G$  ( i.e. if $g = \tau|_V$ for some $\tau \in T$) and if $\alpha \beta \notin E$, then
\begin{align*}
g(\alpha) \wedge g (\beta) \text{  mod  } W &= [\tau(\alpha), \tau(\beta)] \\&= \tau[\alpha, \beta]\\ & = 0,
\end{align*}
 since $\tau$ is an automorphism of $\n$ and $[\alpha, \beta] = 0$ in $\n$. Conversely if $g  \in {\rm GL }(V)$ such that $g(\alpha) \wedge g(\beta) \in W$ for all $\alpha \beta \notin E$, then $\tau = g \oplus \wedge^2 g$ defines a Lie automorphism of $\n$ such that
$\tau(V) = V$ and hence $g \in G$.  Since the condition that an element of $G$ stabilizes $W$ under the natural action on $\bigwedge^2V$ is represented by polynomial equations, we have therefore the following:

\begin{lemma}\label{algebraic}

 $G$ is an algebraic group.

\end{lemma}

In fact,  for any two-step nilpotent Lie algebra $\mathfrak{m} = V \oplus (\bigwedge^2 V)/W$ one can analogously define the subgroup  $G$ of ${\rm GL}(V)$ consisting of the restrictions of automorphisms of $\m$ fixing $V$. A similar argument shows that $G$ is an algebraic group.
However our next lemma holds only for two-step nilpotent Lie algebras arising from graphs.
 We prove that all linear automorphisms of $V$ which can be represented as diagonal matrices with respect to the basis $S$ can be extended to Lie automorphisms of $\n(S, E)$.

\begin{lemma}\label{diagonal}
Let $D_{S}$ denote the subgroup of ${\rm GL}(V)$ consisting of all
 elements   which can be represented
 as diagonal matrices with respect to the basis  $S$ of $V$.
  Then $D_{S} \subseteq G$.
  \end{lemma}

 \proof   We recall that  $G$ is
 a subgroup of  ${\rm GL}(V)$ consisting of those linear automorphisms of $V$ whose induced action on $\bigwedge^2 V$ leaves the subspace $W$ invariant  where $W \subseteq \bigwedge^2V$ spanned by elements $\alpha \wedge \beta$ such that $\alpha \beta \notin E$.
 Let $d \in D_{S}$, say $d(\alpha) = d_{\alpha} \alpha$ for some nonzero $d_{\alpha}$'s
in $k$   and for all $\alpha \in S$.  Then if  $\alpha \beta  \notin E$, we have $d(\alpha) \wedge d(\beta) = d_{\alpha}d_{\beta}(\alpha \wedge \beta) \in W$.  Hence  $d \in G$ and $D_{S} \subseteq G$. \hfill$\square$

 \vspace{.2cm}

\section{Proof of  Theorem~\ref{main}}

Let $\overline{k}$ denote  the  algebraic closure of $k$. First we note that if $F$ is an isomorphism from the Lie algebra $\n = \n(S, E)$ to $\n'= \n(S', E')$, then $ \overline{F} = F  \otimes_k {\rm id}_{\overline{k}}$ ~is an isomorphism of $\overline{k}$-Lie algebras $\overline{\n} = \n \otimes_k \overline{k}$ and $\overline{\n'} = \n' \otimes_k \overline{k}$. To see this, recall that the Lie bracket in $\overline{\n}$ is defined by \[[x \otimes a, y \otimes b] = [x, y] \otimes ab \text{ for all } x, y \in \n \text{ and } a, b \in \overline{k}.\]
 Also $\overline{F}$ is defined by   \[\overline{F} (\sum_i x_i \otimes a_i) = \sum_i F(x_i) \otimes a_i \text{ for all } x_i \in \n \text{ and } a_i \in \overline{k}.\]  Then $\overline{F}$  is a $\overline{k}$-vector space isomorphism which follows from the fact that $F$ is a $k$-vector space isomorphism.  Furthermore, for $x, y \in \n$ and $a, b \in \overline{k}$
 we have
 \begin{align*}
 \overline{F}[x \otimes a, y \otimes b] = \overline{F}([x, y] \otimes ab) &=   (F[x, y] \otimes ab)
\\&=[F(x), F(x)] \otimes ab  \text{ since $F$ is a Lie algebra isomorphism.} \\& =  [F(x) \otimes a , F(x) \otimes b] = [\overline{F}(x \otimes a), \overline{F}(y \otimes b)].
\end{align*}

Now without loss of generality we can assume that the field $k$ is algebraically closed. Indeed if $k$ is not  already algebraically closed, and $F$ is an isomorphism from the Lie algebra $\n = \n(S, E)$ to $\n'= \n(S', E')$,     we can replace  $\n$ by $\overline{\n}$,  $\n'$ by $\overline{\n'}$
 and $F$ by $\overline{F}$. Then the new $F$ is an isomorphism of $\overline{k}$-Lie algebras $\n$ and $\n'$ as discussed above.
  We need to show that the graphs $(S, E)$ and $(S', E')$ defining the Lie algebras $\n$ and $\n'$ are isomorphic.

We recall that $\n = V \oplus (\bigwedge^2V)/W$ where $V$ is the $k$-vector space with basis $S$ and
  $W \subset \bigwedge^2 V$ is as defined in the previous section. Similarly $\n' = V' \oplus (\bigwedge^2V')/W'$  where $V'$ is the $k$-vector space with basis $S'$.

  The proof will be in three  steps.  First we construct a new graph $(S'', E'')$ isomorphic to the graph $(S, E)$. The set of vertices $S''$ will be a
  basis of the vector space $V'$. Then in the next step we look a group $D_{S''}$ analogous to the group $D_S$ considered in Lemma~\ref{diagonal}.
  The theorem will follow from the relation between $D_{S''}$ and $D_{S'}$ to be explained in the third step.

We begin by constructing  a new graph $(S'', E'')$ isomorphic to the graph $(S, E)$ in  the following way.
  Recall that $V'$ denote the $k$-vector space with  $S'$ as basis and let $W'$ denote  the subspace of $\bigwedge^2 V'$ spanned by
 the wedge products $\alpha \wedge \beta$ where $\alpha$ and $ \beta$ are vertices in $S'$  not connected by an edge in $E'$.  Let $\pi: \n' \to V'$ denote
 the canonical linear projection with respect to the decomposition $\n' = V' \oplus (\bigwedge^2V')/W'$.
  We define $S''$ denote the subset of $V'$ given by  $ \{\pi(F(\alpha)) \mid \alpha \in S \}$ and  $E''$  denote   the set of unordered pairs $\pi(F(\alpha)) \pi(F(\beta))$ such that $\alpha
   \beta \in E$, where, recall that $F$ is the isomorphism given between $\n$ and $\n'$.

  It follows that  $S''$ is a basis of $V'$.  Indeed, if $ \sum_{i=1}^{n} a_i \pi(F(\alpha_i)) = 0$ for $a_i$'s $\in \mathbb R$ and $\alpha_i$'s $\in S$,
  then $ F(\sum_{i=1}^{n} a_i \alpha_i) \in [\n', \n']$.
   Since $F$ is a Lie algebra isomorphism, we can see that $ \sum_{i=1}^{n} a_i \alpha_i \in [\n, \n] $. On the other hand each $\alpha_i$ is in $S$ so
   $\sum_{i=1}^{n} a_i \alpha_i \in V \cap [\n, \n]$, so that   $ \sum_{i=1}^{n} a_i \alpha_i = 0$.
By the linear independence of $S$, each $a_i = 0$.
Now the derived algebras $[\n, \n]$ and $[\n', \n']$ have the same dimension since $\n$ and $\n'$ are isomorphic. Hence $|E| = |E'|$ since
$|E| =  \dim [\n, \n]$ and $|E'| = \dim [\n', \n']$. Hence $|S| = |S'|$ because $\n$ is of dimension $|S| + |E|$ and $\n'$ is of dimension $|S'|+|E'|$, and $\n$ and $\n'$ are isomorphic.
Thus  the sets $S'$ and $S''$  have exactly the same number of elements and hence $S''$ is a basis of $V'$ since $S'$ is a basis of $V'$.

Since the graph $(S'', E'')$ as constructed above  is isomorphic to the graph $(S, E)$,  to prove the theorem, it is enough to show that there exists an isomorphism of graphs
$f
: (S', E') \to (S'', E'')$.

Now as in Lemma~\ref{diagonal}, $D_{S'}$ denote the subgroup of ${\rm GL}(V')$ consisting of all
 elements   which can be represented
 as diagonal matrices with respect to the basis  $S'$ of $V'$.
 By Lemma \ref{diagonal}, we have  $D_{S'} \subset G'$.
Similarly let  $D_{S''}$ denote the subgroup of ${\rm GL}(V')$ consisting of all
 elements   which can be represented
 as diagonal matrices with respect to the basis  $S'' =  \{\pi(F(\alpha)) \mid \alpha \in S \}$  of $V'$.

 We claim  that $D_{S''} \subset G'$ as well.   To see this, let $d \in D_{S''}$ and $d(\pi(F(\alpha))
 = d_{\alpha} \pi(F(\alpha)$ for all $\alpha \in S$ and for some nonzero $d_{\alpha}$'s in $k$.
 It suffices to show  that if $\gamma$ and  $\delta$ are vertices in $S'$ not connected by an edge in $E'$, then $d(\gamma) \wedge d(\delta) \in W'$  (i.e. $[d(\gamma), d(\delta)] = 0$ in $\n'$.)
 Now since $\gamma, \delta \in V'$ and $S''$ is a basis for $V'$, we represent  $\displaystyle \gamma=  \sum_{\alpha \in S} a_\alpha \pi(F(\alpha))$ and
 $\displaystyle \delta = \sum_{\beta\in S} b_\beta \pi(F(\beta))$  where  each $a_\alpha$, $b_\beta$ is  in $k$.  Since $[\gamma, \delta] = 0$,  we have
$\displaystyle\left[\pi\left(F\left(\sum_{\alpha \in S} a_\alpha \alpha\right)\right), \pi\left(F\left (\sum_{\beta\in S}  b_\beta \beta\right)\right)\right] = 0.$
This means that \[ \left[F\left(\sum_{\alpha \in S}  a_\alpha \alpha\right), F\left(\sum_{\beta\in S}  b_\beta\beta\right)\right] = 0,\]
because in $\n'$ $[v + w, v' + w'] = 0$ if and only if  $[v, v'] = 0$ for $v, v' \in V'$ and $w, w' \in (\bigwedge^2V')/W'$. Hence
 $\displaystyle \left[\sum_ {\alpha \in S} a_\alpha \alpha, \sum_{\beta\in S} b_\beta \beta\right] = 0$ in $\n$ since $F$ is an isomorphism. Now we define $\sigma$, an element of ${\rm GL}(V)$ by  $\sigma (\alpha) = d_{\alpha} \alpha$ for each $\alpha \in S$.
Then $\sigma \in D_S$ where $D_S$ is the subgroup of ${\rm GL}(V)$ consisting of diagonal automorphisms of $V$ with respect to $S$ and hence $\sigma \in G$ by Lemma \ref{diagonal}. So $\displaystyle \left[\sum_{\alpha\in S} a_\alpha \sigma(\alpha), \sum_{\beta\in S} b_\beta \sigma(\beta)\right] = 0$ in $\n$ which means that  $\displaystyle \left[\sum_{\alpha\in S} a_\alpha d_{\alpha}\alpha, \sum_{\beta\in S} b_\beta d_{\beta} \beta\right] = 0$.
 Since $F$ is an isomorphism, $\displaystyle \left[\sum_{\alpha\in S}  a_\alpha d_{\alpha}F(\alpha), \sum_{\beta\in S} b_\beta d_{\beta} F (\beta)\right] = 0$ in $\n'$ and hence \[\left[\sum_{\alpha\in S} a_\alpha  d_{\alpha}\pi(F(\alpha)), \sum_ {\beta\in S}b_\beta d_{\beta} \pi(F (\beta))\right]= 0.\] Thus $[d(\gamma), d(\delta)] = 0$ in $\n'$ as we wanted to show.

By Lemma~\ref{algebraic}, $G'$ is an algebraic group. Moreover, $D_{S'}$ and $D_{S''}$ are maximal tori in the connected component in Zariski topology
of identity of $G'$ . But maximal tori in a connected algebraic group over an algebraically closed field (such as our $k$)  are conjugate (see \cite{S}, p. 108),
and hence there exists $g \in G'$ such that $$D_{S'} = g\left( D_{S''}\right)g^{-1}.$$

To define an isomorphism of graphs  $f: S' \to S''$, we construct  a special  element of $D_{S'}$. For each $\alpha \in S'$, choose  a nonzero  $d'_{\alpha}\in k$ such that unless the set $\{\alpha, \beta\}$ equals
 the set $\{\gamma, \delta\}$ we have $d'_{\alpha}d'_{\beta} \neq
  d'_{\gamma} d'_{\delta}$.  Since by hypothesis the field $k$ is  algebraically closed, it is in particular infinite, and such  $d_\alpha$'s may be chosen.  Define $d' \in D_{S'}$ be the element whose matrix representation with respect to $S'$ given by
$d'(\alpha) = d'_{\alpha} \alpha$.
 As note above, there exists $d'' \in D_{S''}$ such that $$d'= g d'' g^{-1}.$$
  Since $d'$ and $d''$ are similar, they have the same eigenvalues. Further since $d'$ is diagonal
   with respect to $S'$ and $d''$ is diagonal with respect to $S''$, it follows that these diagonal entries  are the same up to a permutation.
   Hence there exists a bijection $f$ of $S'$
  with $S''$ such that $$d''(f (\alpha)) = d'_{\alpha}
  f(\alpha)  \text { for all }  \alpha \in S'.$$
We claim that $f$ is in fact an isomorphism of the graph $(S', E')$ with the graph $(S'', E'')$. Establishing this claim will conclude the proof of our theorem.

It suffices now to show that  $\alpha \beta \in E'$ if and only if $f(\alpha)f(\beta) \in E''$.
 Now $\alpha \beta \in E' $ if and only if $[\alpha, \beta] \neq 0$ in $\n'$.  Also  $f(\alpha)f(\beta) \in E''$ if and only if $[f(\alpha), f(\beta)] \neq 0$ in $\n'$.
This follows from the definition of $E''$ and using that $F$ is an isomorphism. Indeed $f(\alpha) = \pi(F(\alpha'))$ and $f(\beta) = \pi(F(\beta'))$ for some $\alpha', \beta' \in S$. Hence $f(\alpha)f(\beta) \in E''$ if and only if  $\alpha' \beta' \in E$.
 Since $F$ is an isomorphism, $\alpha' \beta' \in E$ if and only if $[F(\alpha'), F(\beta')] \neq 0$ in $\n'$
   i.e. $[\pi(F(\alpha')), \pi(F(\beta'))] \neq 0$ in $\n'$.

  We will show that for $\alpha, \beta \in S'$,  $[\alpha, \beta] \neq 0$ if and only if $[f(\alpha), f(\beta)] \neq 0$.
   Since the diagonal entries   of the chosen element $d' \in D_{S'}$ (defined as above) are such that the pairwise products of the diagonal entries are
    distinct, $d'_\alpha d'_\beta$ is an eigenvalue of an extended automorphism $d'$ of $\n'$ if and only if $[\alpha, \beta] \neq 0$.
    Similarly $d'_\alpha d'_\beta$ is an eigenvalue of an extended automorphism $d''$ of $\n'$ if and only if $[f(\alpha), f(\beta)] \neq 0$.
       Now  we note that
  \[ \left[g^{-1}(d'(\alpha)), g^{-1}(d'(\beta))] = d'_{\alpha} d'_{\beta} [g^{-1}
(\alpha), g^{-1} (\beta)\right].\]
 On the other hand, we have \[\left[ g^{-1}(d'(\alpha)),
g^{-1}(d('\beta))\right]= d''\left[g^{-1} (\alpha),  g^{-1}(\beta)\right] \text{ since } g^{-1} d' = d'' g^{-1}. \]
Since $g \in G'$, $[g^{-1} (\alpha),  g^{-1}(\beta)]
 \neq 0$ if and only if $[\alpha, \beta] \neq 0$. Hence $d'_{\alpha} d'_{\beta}$ is an eigenvalue of
 $d''$ if and only if $d'_{\alpha} d'_{\beta}$ is an eigenvalue of  $d''$. Thus
  $\alpha \beta \in E'$ if and only if $f(\alpha)f(\beta) \in E''$. Hence $f: (S',E') \to (S'', E'')$ is   an isomorphism of graphs, this complete the proof of Theorem~\ref{main}.

\vspace{0.3cm}

\noindent {\it Acknowledgments:}
  I am thankful to Prof. Jorge Lauret for raising   this question. I am also thankful to Prof. S. G. Dani for his help and Dr.
   Debraj Chakrabarti for help with the final version. I also express my gratitude to the referee for helpful suggestions and pointing out the references \cite{Dr1, DuK1, DuK2, KLNR}.

\end{document}